\input amstex
\input amsppt.sty
\magnification=\magstep1
\hsize=30truecc
\baselineskip=16truept
\vsize=22.2truecm
\nologo
\pageno=1
\topmatter
\TagsOnRight

\def\Z{\Bbb Z}

\def\C{\Bbb C}
\def\al{\alpha}
\def\l{\left}
\def\r{\right}
\def\bg{\bigg}
\def\({\bg(}
\def\){\bg)}
\def\[{\bg[}
\def\]{\bg]}
\def\t{\text}
\def\f{\frac}
\def\em{\emptyset}
\def\se {\subseteq}
\def\sp {\supseteq}
\def\sm{\setminus}

\def\bi{\binom}
\def\eq{\equiv}

\def\ls{\leqslant}
\def\gs{\geqslant}
\def\mo{\t{\rm mod}}

\def\Proof{\noindent{\it Proof}}

\hbox{To appear in Adv. in Appl. Math.}
\bigskip
\title A connection between covers of the integers and unit fractions\endtitle
\rightheadtext{Covers of the integers and unit fractions}
\author Zhi-Wei Sun\endauthor
\affil Department of Mathematics, Nanjing University
\\Nanjing 210093, People's Republic of China
\endaffil
\date Received 8 July 2005; accepted 2 February 2006\enddate

\abstract {For integers $a$ and $n>0$, let $a(n)$ denote the residue class
$\{x\in\Z:\, x\eq a\ (\mo\ n)\}$. Let $A$ be a collection
$\{a_s(n_s)\}_{s=1}^k$ of finitely many residue classes
such that $A$ covers all the integers at least $m$ times
but $\{a_s(n_s)\}_{s=1}^{k-1}$ does not.
We show that if $n_k$ is a period of the covering function
$w_A(x)=|\{1\ls s\ls k:\,x\in a_s(n_s)\}|$ then for any
$r=0,\ldots,n_k-1$ there are at least $m$ integers in the form
$\sum_{s\in I}1/n_s-r/n_k$ with $I\se\{1,\ldots,k-1\}$.}


\endabstract

\thanks {\it E-mail address:} zwsun\@nju.edu.cn
\ \ \ {\it URL:} {\tt http://pweb.nju.edu.cn/zwsun}
\newline\indent Supported by the National Science Fund
for Distinguished Young Scholars (no. 10425103) and a Key
Program of NSF (no. 10331020) in China.
\endthanks
\endtopmatter
\document

\heading{1. Introduction}\endheading

 For an integer $a$ and a positive integer $n$,
we use $a(n)$ to denote the residue class
 $\{x\in\Z:\, x\eq a\ (\mo\ n)\}$. For a finite system
 $$A=\{a_s(n_s)\}_{s=1}^k\tag1$$
 of residue classes, the function $w_A:\Z\to\{0,1,\ldots\}$
 given by
 $$w_A(x)=|\{1\ls s\ls k:\,x\in a_s(n_s)\}|$$
 is called the {\it covering function} of $A$.
 Clearly $w_A(x)$ is periodic modulo the least common multiple
 $N_A$ of the moduli $n_1,\ldots,n_k$, and
 it is easy to verify the following well-known equality:
 $$\f 1{N_A}\sum_{x=0}^{N_A-1}w_A(x)=\sum_{s=1}^k\f1{n_s}.$$
 As in [7] we call $m(A)=\min_{x\in\Z}w_A(x)$
 the {\it covering multiplicity} of $A$.
 For example,
 $$B=\{0(2),\ 0(3),\ 1(4),\ 5(6),\ 7(12)\}$$
 has covering multiplicity $m(B)=1$, because the covering function is periodic modulo $N_B=12$,
 and $$w_B(x)=\cases1&\t{if}\ x\in\{1,2,3,4,7,8,10,11\},
 \\2&\t{if}\ x\in\{0,5,6,9\}.\endcases$$

 Let $m$ be a positive integer.
 If $w_A(x)\gs m$ for all $x\in\Z$, then we call $A$ an {\it $m$-cover} of the integers,
 and in this case we have the well-known inequality $\sum_{s=1}^k1/n_s\gs m$.
 (The term ``1-cover" is usually replaced by the word ``cover".)
 If $A$ is an $m$-cover of the integers but $A_t=\{a_s(n_s)\}_{s\in[1,k]\sm\{t\}}$ is not
 (where $[a,b]=\{x\in\Z:\,a\ls x\ls b\}$ for $a,b\in\Z$),
 then we say that
 $A$ forms an $m$-cover of the integers with $a_t(n_t)$ {\it irredundant}.
 (For example, $\{0(2),1(2),2(3)\}$ is a cover of the integers in which
 $2(3)$ is redundant while $0(2)$ and $1(2)$ are irredundant.)
 If $w_A(x)=m$ for all $x\in\Z$,
 then $A$ is said to be an {\it exact $m$-cover} of the integers,
 and in this case we have the equality $\sum_{s=1}^k1/n_s=m$.
 By a graph-theoretic argument, M. Z. Zhang [12] showed that if $m>1$ then
 there are infinitely many exact $m$-covers of the integers
 each of which cannot be split into an $n$-cover
 and an $(m-n)$-cover of the integers with $0<n<m$.

 Covers of the integers by residue classes
 were first introduced by P. Erd\H os (cf. [1]) in the 1930s, who observed that
 the system $B$ mentioned above
 is a cover of the integers with distinct moduli.
 The topic of covers of the integers has been an active one
 in combinatorial number theory (cf. [3, 4]),
 and many surprising applications have been found
 (see, e.g., [1, 2, 9, 11]).
 The so-called $m$-covers and exact $m$-covers of the integers
 were systematically studied
 by the author in the 1990s.

 Concerning the cover $B$ given above one can easily check that
 $$\align&\l\{\sum_{n\in S}\f1n:\, S\se\{2,3,4,6,12\}\r\}
 \\=&\l\{0,\f1{12},
 \ldots,\f{11}{12}\r\}\bigcup\l\{1+\f r{12}:\ r=0,1,2,3,4\r\}.
 \endalign$$
 This suggests that for a general $m$-cover (1) of the integers
 we should investigate the set
 $\{\sum_{s\in I}1/n_s:\, I\se[1,k]\}$.

 In this paper we establish the following new connection
 between covers of the integers and unit fractions.

 \proclaim{Theorem 1} Let $A=\{a_s(n_s)\}_{s=1}^k$ be an $m$-cover of the integers with
 the residue class $a_k(n_k)$ irredundant.
 If the covering function $w_A(x)$ is periodic modulo $n_k$,
 then for any $r=0,\ldots,n_k-1$ we have
 $$\bg|\bg\{\bg\lfloor\sum_{s\in I}\f1{n_s}\bg\rfloor:\
I\se[1,k-1]\ \t{and}\ \bg\{\sum_{s\in I}\f1{n_s}\bg\}=\f r{n_k}\bg\}\bg|\gs m,\tag2$$
where $\lfloor\al\rfloor$ and $\{\al\}$ denote the integral part and
the fractional part of a real number $\al$
respectively.
\endproclaim

 Note that $n_k$ in Theorem 1 needn't be the largest modulus
among $n_1,\ldots,n_k$.
In the case $m=1$ and $n_k=N_A$, Theorem 1 is an easy consequence
of [5, Theorem 1] as observed by the author's twin brother Z. H. Sun.
When $w_A(x)=m$ for all $x\in\Z$, the author [6] even proved the following stronger result:
$$\bg|\bg\{I\se[1,k-1]:\, \sum_{s\in I}\f 1{n_s}=\f a{n_k}\bg\}\bg|
\gs\bi{m-1}{\lfloor a/n_k\rfloor}
\quad\ \t{for all}\ a=0,1,\ldots.$$
Given an $m$-cover $\{a_s(n_s)\}_{s=1}^k$ of the integers with $a_k(n_k)$
 irredundant, by refining a result in [7] the author can show that
 there exists a real number $0\ls\al<1$
 such that (2) with $r/n_k$ replaced by $(\al+r)/n_k$ holds
 for every $r=0,\ldots,n_k-1$.

 Here we mention two local-global results related to Theorem 1.

(a) (Z. W. Sun [5]) $\{a_s(n_s)\}_{s=1}^k$ forms an $m$-cover of the integers
if it covers $|\{\{\sum_{s\in I}1/n_s\}:\,I\se[1,k]\}|$
consecutive integers at least $m$ times.

(b) (Z. W. Sun [10]) $\{a_s(n_s)\}_{s=1}^k$ is an exact $m$-cover of the integers if it covers
$|\bigcup_{s=1}^k\{r/n_s:\,r\in[0,n_s-1]\}|$ consecutive integers exactly $m$ times.

\proclaim{Corollary 1} Suppose that the covering function
of $A=\{a_s(n_s)\}_{s=1}^k$ has a positive integer period $n_0$.
If there is a unique $a_0\in[0,n_0-1]$ such that
$w_A(a_0)=m(A)$, then for any $D\se\Z$ with $|D|=m(A)$ we have
$$\bg\{\bg\{\sum_{s\in I}\f1{n_s}\bg\}:\, I\se[1,k]\ \t{and}
\ \bg\lfloor\sum_{s\in I}\f1{n_s}\bg\rfloor\not\in D\bg\}
\sp\bg\{\f r{n_0}:\, r\in[0,n_0-1]\bg\}.$$
\endproclaim
\Proof. Let $m=m(A)+1$. Clearly $A'=\{a_s(n_s)\}_{s=0}^k$
forms an $m$-cover of the integers with $a_0(n_0)$ irredundant.
As $w_{A'}(x)-w_A(x)$ is the characteristic function of $a_0(n_0)$,
$w_{A'}(x)$ is also periodic mod $n_0$. Applying Theorem 1 we immediately get
the desired result. \qed

\medskip
We will prove Theorem 1 in Section 3 with help from some lemmas given in the next section.

\heading{2. Several Lemmas}\endheading

\proclaim{Lemma 1} Let $(1)$ be a finite system of residue classes with $m(A)=m$,
and let $m_1,\ldots,m_k$ be any integers.
If $f(x_1,\ldots,x_k)$ is a polynomial with coefficients
in the complex field $\C$ and $\deg f\ls m$,
then for any $z\in\Z$ we have
$$\aligned&\sum_{I\se[1,k]}(-1)^{|I|}f([\![1\in I]\!],\ldots,[\![k\in I]\!])
e^{2\pi i\sum_{s\in I}(a_s-z)m_s/n_s}
\\&\quad=(-1)^kc(I_z)\prod_{s\in[1,k]\sm I_z}\l(e^{2\pi i(a_s-z)m_s/n_s}-1\r),
\endaligned\tag3$$
where
$[\![s\in I]\!]$ takes $1$ or $0$ according as $s\in I$ or not,
$I_z=\{1\ls s\ls k:\, z\in a_s(n_s)\}$,
and $c(I_z)=[\prod_{s\in I_z}x_s]f(x_1,\ldots,x_k)$
is the coefficient of the monomial $\prod_{s\in I_z}x_s$ in $f(x_1,\ldots,x_k)$.
\endproclaim
\Proof. Write $f(x_1,\ldots,x_k)
=\sum_{j_1,\ldots,j_k\gs0}c_{j_1,\ldots,j_k}x_1^{j_1}\cdots x_k^{j_k}$. Observe that
$$\align &\sum_{I\se[1,k]}(-1)^{|I|}f([\![1\in I]\!],\ldots,[\![k\in I]\!])
e^{2\pi i\sum_{s\in I}(a_s-z)m_s/n_s}
\\=&\sum\Sb j_1,\ldots,j_k\gs0\\j_1+\cdots+j_k\ls m\endSb
c_{j_1,\ldots,j_k}\sum_{I\se[1,k]}\(\prod_{s=1}^k[\![s\in I]\!]^{j_s}\times
(-1)^{|I|}e^{2\pi i\sum_{s\in I}(a_s-z)m_s/n_s}\)
\\=&\sum\Sb j_1,\ldots,j_k\gs0\\j_1+\cdots+j_k\ls m\endSb
c_{j_1,\ldots,j_k}\sum_{J(j_1,\ldots,j_k)\se I\se[1,k]}
(-1)^{|I|}e^{2\pi i\sum_{s\in I}(a_s-z)m_s/n_s},
\endalign$$
where $J(j_1,\ldots,j_k)=\{1\ls s\ls k:\,j_s\not=0\}$.

Let $z$ be any integer. If $I_z\not\se J(j_1,\ldots,j_k)$, then
$$\sum_{J(j_1,\ldots,j_k)\se I\se[1,k]}
(-1)^{|I|}e^{2\pi i\sum_{s\in I}(a_s-z)m_s/n_s}=0$$
since
$$\align&\sum_{I\se[1,k]\sm J(j_1,\ldots,j_k)}(-1)^{|I|}e^{2\pi i\sum_{s\in I}(a_s-z)m_s/n_s}
\\=&\prod_{s\in [1,k]\sm J(j_1,\ldots,j_k)}\l(1-e^{2\pi i(a_s-z)m_s/n_s}\r)=0.
\endalign$$
If $j_1,\ldots,j_k$ are nonnegative integers with
$j_1+\cdots+j_k\ls m$ and $I_z\se J(j_1,\ldots,j_k)$, then
$$m=m(A)\ls|I_z|\ls|J(j_1,\ldots,j_k)|\ls j_1+\cdots+j_k\ls m;$$
hence $I_z=J(j_1,\ldots,j_k)$ and $j_s=1$ for all $s\in I_z$.

Combining the above we find that the left-hand side of (3) coincides with
$$\align &c(I_z)\sum_{I_z\se I\se[1,k]}(-1)^{|I|}e^{2\pi i\sum_{s\in I}(a_s-z)m_s/n_s}
\\=&c(I_z)(-1)^{|I_z|}e^{2\pi i\sum_{s\in I_z}(a_s-z)m_s/n_s}
\prod_{s\in[1,k]\sm I_z}\l(1-e^{2\pi i(a_s-z)m_s/n_s}\r)
\\=&(-1)^kc(I_z)\prod_{s\in [1,k]\sm I_z}\l(e^{2\pi i(a_s-z)m_s/n_s}-1\r).
\endalign$$
This proves the desired (3).
 \qed

\proclaim{Lemma 2} Let $(1)$ be an $m$-cover of the integers with $a_k(n_k)$ irredundant,
and let $m_1,\ldots,m_{k-1}$ be positive integers.
Then, for any $0\ls\al<1$ we have
$C_0(\al)=\cdots=C_{n_k-1}(\al)$, where $C_r(\al)\ ($with $r\in[0,n_k-1])$ denotes the sum
$$\sum\Sb I\se[1,k-1]\\\{\sum_{s\in I}m_s/n_s\}=(\al+r)/n_k\endSb(-1)^{|I|}
\bi{\lfloor\sum_{s\in I}m_s/n_s\rfloor}{m-1}e^{2\pi i\sum_{s\in I}(a_s-a_k)m_s/n_s}.$$
\endproclaim
\Proof. This follows from [7, Lemma 2]. \qed

\proclaim{Lemma 3} Let $(1)$ be an $m$-cover of the integers with $a_k(n_k)$ irredundant.
Suppose that $n_k$ is a period of the covering function $w_A(x)$. Then, for any
$z\in a_k(n_k)$ we have
$$\prod_{s\in[1,k]\sm I_z}\l(1-e^{2\pi i(a_s-z)/n_s}\r)
=\prod_{s\in I_z}n_s\times\prod_{t=1}^{n_k}\l(1-e^{2\pi i(t-a_k)/n_k}\r)^{w_A(t)-m}$$
where $I_z=\{1\ls s\ls k:\, z\in a_s(n_s)\}$.
\endproclaim
\Proof. Since $a_k(n_k)$ is irredundant, we have $w_A(z_0)=m$ for some $z_0\in a_k(n_k)$.
As the covering function of $A$ is periodic mod $n_k$, $|I_z|=w_A(z)=m$ for all $z\in a_k(n_k)$.

Now fix $z\in a_k(n_k)$.
Since $w_A(x)$ is periodic modulo $n_k$, by [8, Lemma 2.1] we have the identity
$$\prod_{s=1}^k\l(1-y^{N/n_s}e^{2\pi i a_s/n_s}\r)
=\prod_{t=1}^{n_k}\l(1-y^{N/n_k}e^{2\pi it/n_k}\r)^{w_A(t)},$$
where $N=N_A$ is the least common multiple of $n_1,\ldots,n_k$.
Putting $y=r^{1/N}e^{-2\pi iz/N}$ where $r\gs0$, we then get that
$$\prod_{s=1}^k\l(1-r^{1/n_s}e^{2\pi i(a_s-z)/n_s}\r)
=\prod_{t=1}^{n_k}\l(1-r^{1/n_k}e^{2\pi i(t-z)/n_k}\r)^{w_A(t)}.$$
Therefore
$$\align&\prod_{s\in[1,k]\sm I_z}\l(1-e^{2\pi i(a_s-z)/n_s}\r)
\\=&\lim_{r\to1}\prod_{s\in[1,k]\sm I_z}\l(1-r^{1/n_s}e^{2\pi i(a_s-z)/n_s}\r)
\\=&\lim_{r\to1}\f{\prod_{t=1}^{n_k}(1-r^{1/n_k}e^{2\pi i(t-a_k)/n_k})^{w_A(t)}}
{\prod_{s\in I_z}(1-r^{1/n_s}e^{2\pi i(a_s-z)/n_s})}
\\=&\lim_{r\to1}\prod_{s\in I_z}\f{1-r}{1-r^{1/n_s}}
\times\lim_{r\to 1}\f{\prod_{t=1}^{n_k}(1-r^{1/n_k}e^{2\pi i(t-a_k)/n_k})^{w_A(t)}}{(1-r)^m}
\\=&\prod_{s\in I_z}n_s\times\lim_{r\to 1}\prod_{t=1}^{n_k}(1-r^{1/n_k}
e^{2\pi i(t-a_k)/n_k})^{w_A(t)-m},
\endalign$$
and hence the desired result follows. \qed

\heading{3. Proof of Theorem 1}\endheading

 In the case $k=1$, we must have  $m=1$ and $n_k=1$;
hence the required result is trivial. Below we assume that $k>1$.

Let $r_0\in[0,n_k-1]$ and
$D=\{d_n+r_0/n_k:\, n\in[1,m-1]\}$, where $d_1,\ldots,d_{m-1}$ are
$m-1$ distinct nonnegative integers. (If $m=1$ then we set $D=\em$.)
We want to show that there exists an $I\se[1,k-1]$ such that
$\{\sum_{s\in I}1/n_s\}=r_0/n_k$ and $\sum_{s\in I}1/n_s\not\in D$.

Define
$$f(x_1,\ldots,x_{k-1})=\prod_{d\in D}\l(\f{x_1}{n_1}+\cdots+\f{x_{k-1}}{n_{k-1}}-d\r).$$
(An empty product is regarded as 1.) Then $\deg f=|D|=m-1$. For
any $z\in a_k(n_k)$, the set $I_z=\{1\ls s\ls k:\,z\in a_s(n_s)\}$
has cardinality $m$ since $a_k(n_k)$ is irredundant and $w_A(x)$
is periodic mod $n_k$. Observe that the coefficient
$$c_z=\[\prod_{s\in I_z\sm\{k\}}x_s\]f(x_1,\ldots,x_{k-1})
=\[\prod_{s\in I_z\sm\{k\}}x_s\]\(\sum_{s=1}^{k-1}\f{x_s}{n_s}\)^{m-1}$$
coincides with $(m-1)!/\prod_{s\in I_z\sm\{k\}}n_s$ by the multinomial theorem.
For $I\se[1,k-1]$ we set
$$v(I)=f([\![1\in I]\!],\ldots,[\![k-1\in I]\!]).$$
As $|\{1\ls s\ls k-1:\, x\in a_s(n_s)\}|\gs \deg f$ for all $x\in\Z$,
in view of Lemmas 1 and 3 we have
$$\align &\sum_{I\se[1,k-1]}(-1)^{|I|}v(I)e^{2\pi i\sum_{s\in I}(a_s-z)/n_s}
\\=&(-1)^{k-1}c_z\prod_{s\in[1,k-1]\sm I_z}\l(e^{2\pi i(a_s-z)/n_s}-1\r)
\\=&\f{(-1)^{m-1}(m-1)!}{\prod_{s\in I_z\sm\{k\}}n_s}{\prod_{s\in I_z}n_s}
\times\prod_{t=1}^{n_k}\l(1-e^{2\pi i(t-a_k)/n_k}\r)^{w_A(t)-m}=C,
\endalign$$
where $C$ is a nonzero constant not depending on $z\in a_k(n_k)$.

By the above,
$$\align N_AC=&\sum_{x=0}^{N_A-1}\sum_{I\se[1,k-1]}(-1)^{|I|}v(I)
e^{2\pi i\sum_{s\in I}(a_s-a_k-n_kx)/n_s}
\\=&\sum_{I\se[1,k-1]}(-1)^{|I|}v(I)e^{2\pi i\sum_{s\in I}(a_s-a_k)/n_s}
\sum_{x=0}^{N_A-1}e^{-2\pi ix\sum_{s\in I}n_k/n_s}
\endalign$$
and hence
$$C=\sum\Sb I\se[1,k-1]\\n_k\sum_{s\in I}1/n_s\in\Z\endSb(-1)^{|I|}v(I)
e^{2\pi i\sum_{s\in I}(a_s-a_k)/n_s}
=\sum_{r=0}^{n_k-1}C_r,$$
where
$$C_r=\sum\Sb I\se[1,k-1]\\\{\sum_{s\in I}1/n_s\}=r/n_k\endSb
(-1)^{|I|}\prod_{d\in D}\(\sum_{s\in I}\f1{n_s}-d\)
e^{2\pi i\sum_{s\in I}(a_s-a_k)/n_s}.$$

 Let $r\in[0,n_k-1]$. Write
 $$P_r(x)=\prod_{d\in D}\l(x+\f r{n_k}-d\r)=\sum_{n=0}^{m-1}c_{n,r}\bi xn$$
 where $c_{n,r}\in\C$. By comparing the leading coefficients, we find that
 $c_{m-1,r}=(m-1)!$. Observe that
 $$\align C_r=&\sum\Sb I\se[1,k-1]\\\{\sum_{s\in I}1/n_s\}=r/n_k\endSb
 (-1)^{|I|}P_r\(\bg\lfloor\sum_{s\in I}\f1{n_s}\bg\rfloor\)e^{2\pi i\sum_{s\in I}(a_s-a_k)/n_s}
 \\=&\sum_{n=0}^{m-1}c_{n,r}\sum\Sb I\se[1,k-1]\\\{\sum_{s\in I}1/n_s\}=r/n_k\endSb
 (-1)^{|I|}\bi{\lfloor\sum_{s\in I}1/n_s\rfloor}ne^{2\pi i\sum_{s\in I}(a_s-a_k)/n_s}
 \\=&c_{m-1,r}\sum\Sb I\se[1,k-1]\\\{\sum_{s\in I}1/n_s\}=r/n_k\endSb
 (-1)^{|I|}\bi{\lfloor\sum_{s\in I}1/n_s\rfloor}{m-1}e^{2\pi i\sum_{s\in I}(a_s-a_k)/n_s};
 \endalign$$
 in taking the last step we note that if $0\ls n<m-1$ then
$$\sum\Sb I\se[1,k-1]\\\{\sum_{s\in I}1/n_s\}=r/n_k\endSb
 (-1)^{|I|}\bi{\lfloor\sum_{s\in I}1/n_s\rfloor}{n}e^{2\pi i\sum_{s\in I}a_s/n_s}=0$$
 by [5, Theorem 1] (since $\{a_s(n_s)\}_{s=1}^{k-1}$ is an $(m-1)$-cover of the integers).
 By Lemma 2 and the above,
 $$C_r=(m-1)!\sum\Sb I\se[1,k-1]\\\{\sum_{s\in I}1/n_s\}=0\endSb
 (-1)^{|I|}\bi{\lfloor\sum_{s\in I}1/n_s\rfloor}{m-1}e^{2\pi i\sum_{s\in I}(a_s-a_k)/n_s}$$
 does not depend on $r\in[0,n_k-1]$.

  Combining the above we obtain that
 $$n_kC_{r_0}=\sum_{r=0}^{n_k-1}C_r=C\not=0.$$
 So there is an $I\se[1,k-1]$ for which $\{\sum_{s\in I}1/n_s\}=r_0/n_k$,
$\sum_{s\in I}1/n_s\not\in D$ and hence
$\lfloor\sum_{s\in I}1/n_s\rfloor\not\in\{d_n:\,n\in[1,m-1]\}$.
This concludes our proof.

\enddocument